\documentclass[12pt]{amsart}
\usepackage{amsmath}
\usepackage{amsfonts,amssymb,amsthm}
\newtheorem{theorem}{Theorem}[section]

\newtheorem{remark}[theorem]{Remark}
\newtheorem{definition}[theorem]{Definition}

\begin{document}
\title{A multidimensional Law of Sines}
\author{Igor Rivin}
\address{Mathematics Department, Temple University, 1805 N Broad St\\
  Philadelphia, PA 19122}
\address{Mathematics Department, Princeton University \\
Fine Hall,   Washington Rd \\
Princeton, NJ 08544}
\email{rivin@math.temple.edu}
\thanks{The author was partially supported by the Department of
  Mathematical Sciences of the National Science Foundation. The author would also like to thank M.~Pashkevich and
  the other participants of the geometry seminar at the Novosibirsk State University for careful reading of a previous
  version of this paper.}
\begin{abstract}
We give a linear-algebraic proof of the \emph{law of sines}, which
also allows us to extend this theorem to simplices in
$\mathbb{E}^n,$ as Theorem \ref{multisine}.
\end{abstract}
\maketitle Let $ABC$ be a triangle in the euclidean plane. The
classical Law of Sines states that
\begin{equation}
\label{thos}
\frac{|AB|}{\sin \gamma} = \frac{|AC|}{\sin \beta} = \frac{|BC|}{\sin
  \alpha}.
\end{equation}
There are many proofs of this fact, some of them found in high
school geometry textbook, but in this note we will derive this
result as a special case of a theorem about simplices in
$\mathbb{E}^n$ (Theorem \ref{multisine}), which, despite the fact
that the classical Law of Sines has been known for at least two
thousand years, seems to have not been noticed to date.

To begin, let $\Delta$ be a simplex in $\mathbb{E}^n.$ This simplex
will have $n+1$ faces $f_1, \dots, f_{n+1},$ and we will denote the
($n-1$-dimensional) area of $f_i$ by $A_i$ and the \emph{outward unit
  normal} to $f_i$ by $\mathbf{f}_i.$
The following fact is fundamental:
\begin{theorem}
\label{divergence}
$$
\sum_{i=1}^{n+1} A_i \mathbf{f}_i = 0.
$$
\end{theorem}
\begin{proof}
A vector $v$ in $\mathbb{E}^n$ is $0$ if and only if its scalar product
with any other vector $w$ is $0.$ Now let
$$
v = \sum_{i=1}^{n+1} A_i \mathbf{f}_i.
$$
Then,
$$\langle v, w\rangle = \sum{i=1}^{n+1} A_i \langle \mathbf{f}_i, w
\rangle.$$
Each summand $A_i \langle \mathbf{f}_i, w \rangle$ is simply the
\emph{signed} area of the projection of the face $f_i$ onto the plane
$P_w$ through the origin normal to $w.$ Since almost every point of $P_w$
either does not lie in the image of $\Delta$ under the orthogonal
projection in the direction $w$ or is covered twice, with opposing
signs, the result follows.
\end{proof}
\begin{remark}
The above is a special case of Stokes' formula.
\end{remark}
We now construct the matrix of columns
$$\mathcal{F}_{\Delta} = \left(\mathbf{f}_1, \dots, \mathbf{f_{n+1}}\right),$$
and then the \emph{Gram matrix}
$$G_{\Delta} = \mathcal{F}_{\Delta}^t \mathcal{F}_{\Delta}.$$
(from here on we will drop the subscript $\Delta.$) The $ij$-th entry
of $G$ is simply $\langle \mathbf{f}_i, \mathbf{f}_j \rangle,$ which
is the cosine of the exterior dihedral angle between $f_i$ and $f_j.$

If $\Delta$ is a non-degenerate simplex, then the rank of
$\mathcal{F}$ is equal to $n$ (the one linear relation between the
$\mathbf{f}_i$ is given by the Theorem \ref{divergence}) and $G$ is a
symmetric matrix with $n$ positive eigenvalues, and one $0$
eigenvalue.
\begin{theorem}
The null-space of $G$ is spanned by the vector $\mathbf{a} = \left(A_1, \dots, A_{n+1}\right).$
\end{theorem}
\begin{proof}
$$G \mathbf{a} = \mathcal{F}^t \mathcal{F} \mathbf{a} =
\mathcal{F}^t \left(\left\langle \sum_{i=1}^{n+1} A_i \mathbf{f_i}, e_1\right\rangle,
\dots,
\left\langle \sum_{i=1}^{n+1} A_i \mathbf{f_i}, e_n\right\rangle\right)
= \mathbf{0}.
$$
\end{proof}
At this point we need some linear algebra:
\section{Some facts about matrices}
\begin{definition}
Let $M$ be a matrix. The \emph{adjugate} $\widehat{M}$ of $M$ is the matrix of
cofactors of $M.$ That is, $\widehat{M}_{ij} = (-1)^{i+j} \det M^{ij},$
where $M^{ij}$ is $M$ with the $i$-th row and $j$-th column removed.
\end{definition}

The reason for this definition is
\begin{theorem}[Cramer's rule]
For any $n\times n$ matrix $M$ (over any commutative ring)
\begin{equation*}
M \widehat{M} = \widehat{M} M = (\det M) I(n),
\end{equation*}
where $I(n)$ is the $n \times n$ identity matrix.
\end{theorem}

We also need
\begin{definition}
The \emph{outer product} of column vectors $v = (v_1, \dots, v_n)$
and $w = (w_1, \dots, w_n)$ is the matrix $v w^t.$
\end{definition}
Consider an arbitrary vector $x = (x_1, \dots, x_n).$ We see that
\begin{equation}
\left[(v w^t) x\right]_k = \sum_{i=1}^n v_k w_i x_i =
v_k \langle w, x \rangle,
\end{equation}
so that
\begin{equation}
\label{outie}
(v w^t) x = \langle w, x \rangle v.
\end{equation}
We see that $v w^t$ is a multiple of the projection operator
onto the subspace spanned by $v.$ In particular, in the case when $\|
v \| = 1,$ the operator $v v^t$ is the orthogonal projection
operator onto the subspace spanned by $v.$ Since $v w^t $ is a
rank $1$ operators all but one of its eigenvalues are equal to $0.$ The one
(potentially) nonzero eigenvalue equals $\langle v, w\rangle.$

We now show:

\begin{theorem}
\label{nulladj}
Suppose that $M$ has nullity $1$, and the null space of $M$ is
spanned by the vector $v,$ while the null space of $M^t$ is spanned by
the vector $w.$ Then
$$
\widehat{M}= c v w^t,
$$
\end{theorem}
\begin{proof}
Since $M$ is singular, we know that $\det M = 0,$ and so every column
of $\widehat M$ is in the null-space of $M.$ so, letting
$\mathbf{m}_i$ denote the $i$th column of $\widehat{M},$ we see that
\begin{equation*}
\mathbf{m}_i = d_i v.
\end{equation*}
However, $\widehat{M^t} = \left(\widehat{M}\right)^t$  so performing
the computation on transposes we see that
\begin{equation*}
\mathbf{m}_i^t = e_i w.
\end{equation*}
We see that
\begin{equation*}
\widehat{M}_{ij} = d_j v_i = e_i w_j.
\end{equation*}
Writing $d_i = g_i w_i,$ and $e_j = h_j v_j,$ we see that, for
every pair $i, j,$ $g_i w_i v_j = h_j w_i v_j.$ Hence $g_i = h_j =
c,$ and the conclusion follows.
\end{proof}

\begin{theorem}
The constant $c$ in the statement of the last theorem equals the
product of the nonzero eigenvalues of $M$ divided by the inner product
of $v$ and $w.$
\end{theorem}
\begin{proof}
By considering the characteristic polynomial of $M$ we see that the
product of the nonzero eigenvalues of $M$ equals the sum of the
principal $n-1$  minors. On the other hand, the principal
minors of $M$ equal the diagonal elements of $\widehat{M},$ so
\begin{equation*}
c \sum_{i=1}^n w_i v_i = \prod_{i=1}^{n-1} \lambda_j.
\end{equation*}
\end{proof}
\begin{remark}
\label{regdet} By the discussion following Eq. (\ref{outie}), the
product of nonzero eigenvalues of $M$ equals
$$\frac{\det \left(M + w \otimes v\right)}{\langle v, w\rangle}.$$
\end{remark}
\section{Back to simplices}
Let us now return to the Gram matrix $G$ of a simplex $\Delta.$
The results in the preceeding section, combined with Theorem
\ref{divergence} immediately imply:
\begin{theorem}[Multidimensional theorem of sines]
\label{multisine}
For any $1\leq i, j, k, l \leq n+1$
\begin{equation*}
\frac{A_iA_j}{A_kA_l} = \frac{\widehat{G}_{ij}}{\widehat{G}_{kl}},
\end{equation*}
where $A_i, A_j, A_k, A_l$ refer to the areas of the corresponding
faces of $\Delta$ and $\widehat{G}_{ij}$ is the $ij$-th minor of the
Gram matrix of $\Delta.$
\end{theorem}
\begin{proof}
Immediate from Theorem \ref{nulladj} together with Theorem \ref{divergence}.
\end{proof}
\section{Examples} In two dimensions, the Gram matrix of a triangle
$ABC$ is
$$
G_{ABC} = \left( \begin{array}{rrr}
1 & -\cos \gamma & -\cos \beta\\
-\cos \gamma & 1 & -\cos \alpha\\
-\cos \beta & - \cos \alpha & 1
\end{array} \right).
$$

Thus, $\widehat{G}_{11} = \sin^2 \alpha,$ while $\widehat{G}_{22} =
\sin^2 \beta,$ so Theorem \ref{multisine} implies that
$$\frac{|BC|^2}{|AC|^2} = \frac{\sin^2 \alpha}{\sin^2 \beta}.$$
The sign indeterminacy due to the squares is illusory, since all the
quantities involved are \emph{a priori} positive.

Note further that $\widehat{G}_{12} = - \cos \alpha \cos \beta -
\cos \gamma,$ So Theorem \ref{multisine} implies that
$$\frac{|BC|}{|AC|} = \frac{\sin^2 \alpha}{-\cos \alpha \cos \beta -
  \cos \gamma},$$
which, together with the Theorem of Sines proves either the addition
  formula for cosine (if we assume that the sum of the angles of a
  triangle is $\pi$) or that the sum of the angles of a triangle is
  $\pi$ (if we assume the addition formula for cosine).

In three dimensions, for an arbitrary (nondegenerate) tetrahedron
$\Delta,$
$$
G_{\Delta} = \left( \begin{array}{rrrr}
1 & -\cos \alpha_{12} & -\cos \alpha_{13} & -\cos \alpha_{14}\\
-\cos \alpha_{12} & 1 & -\cos \alpha_{23} & -\cos \alpha_{24}\\
-\cos \alpha_{13} & -\cos \alpha_{23} & 1 & -\cos \alpha_{34}\\
-\cos \alpha_{14} & -\cos \alpha_{24} & -\cos \alpha_{34} & 1
\end{array} \right).
$$
A quick computation shows that
\begin{equation}
\label{sineq}
 \frac{A_4^2}{A_3^2} = \frac{1- \cos^2 \alpha_{12}  -
\cos^2 \alpha_{13} - \cos^2
  \alpha_{23} - 2\cos \alpha_{12} \cos \alpha_{23} \cos \alpha_{13}}{1-
  \cos^2 \alpha_{12}  - \cos^2 \alpha_{14} - \cos^2
  \alpha_{24} - 2\cos \alpha_{12} \cos \alpha_{24} \cos \alpha_{14}},
\end{equation}
while using Theorem \ref{multisine} with non-principal minors
shows relationships analogous to the angle-sum relationship for a
triangle.

Marina Pashkevich has pointed out that Eq. (\ref{sineq}) can be
considerably simplified (using the spherical Heron's formula), as
follows:
\begin{equation}
\label{pashkevich}
\frac{A_4}{A_3} = \frac{\sin(S_4/2) \cos(\alpha_{13}/2)
\cos(\alpha_23/2)}{\sin(S_3/2)
\cos(\alpha_{14}/2)\cos(\alpha_{24}/2)},
\end{equation}
where $S_i$ is the (spherical) area of the link of the $i$-th
vertex. M.~Pashkevich has also succeeded in extending this formula
(in $3$ dimensions) to hyperbolic and spherical simplices.
\bibliographystyle{alpha}
\end{document}